\documentclass{aastex62}
\usepackage{nicefrac}
\shorttitle{Envelope Polyhedra}
\shortauthors{Gott}

\begin{document}

\title{Envelope Polyhedra}

\correspondingauthor{J. Richard Gott, III}
\email{jrg@astro.princeton.edu}

\author[0000-0002-0786-7307]{J. Richard Gott, III}
\affil{Department of Astrophysical Sciences \\
Princeton University \\
Princeton, NJ 08544, USA}



\begin{abstract}

This paper presents an additional class of regular polyhedra---envelope
polyhedra---made of regular polygons, where the arrangement of polygons
(creating a single surface) around each vertex is identical; but
dihedral angles between faces need not be identical, and some of the
dihedral angles are $0^\circ$ (i.e., some polygons are placed back to
back). For example, {\em squares--6 around a point} \{4,6\} is produced by
deleting the triangles from the rhombicuboctahedron, creating a hollow
polyhedron of genus 7 with triangular holes connecting 18 interior and
18 exterior square faces.  An empty cube missing its top and bottom
faces becomes an envelope polyhedron, {\em squares--4 around a point} \{4,4\}
with a toroidal topology. This definition leads to many interesting
finite and infinite multiply connected regular polygon networks,
including one infinite network with {\em squares--14 around a point} \{4,14\}
and another with {\em triangles--18 around a point} \{3,18\}.  These are
introduced just over 50 years after my related paper on infinite
spongelike pseudopolyhedra in American Mathematical Monthly (Gott,
1967).

\end{abstract}



\section{Introduction---Pseudopolyhedra} \label{sec:intro}

My work on envelope polyhedra grows directly out of my earlier work on
pseudopolyhedra, which I will describe first. This was my high school
science fair project which won 1st Place in mathematics at the (May,
1965) National Science Fair International (now the Intel International
Science and Engineering Fair). A picture of this project appears in my
book {\it The Cosmic Web} (2016), along with a description.  These were
infinite spongelike polyhedra whose polygons were all regular, whose
vertices were congruent, and two polygons always shared only one
edge. All had a sum of polygon angles around a vertex $> 360^\circ$, and
corresponded to surfaces with negative curvature.

Positively curved surfaces like the sphere can be approximated by
regular polyhedra where the sum of angles at a vertex is $< 360^\circ$.  The
cube is a very rough approximation to a sphere.  A cube is made up of
squares meeting 3 around a point. At the corner of a cube, three
square faces meet at a vertex, and each square has a $90^\circ$ angle
at its corner, making the sum of the angles around a point at the
corner $3 \times 90^\circ$, or $270^\circ$.  This is $90^\circ$ less
than $360^\circ$.  All 5 regular Platonic polyhedra: tetrahedron
({\em triangles--3 around a point}), octahedron ({\em triangles--4 around a
point}), icosahedron ({\em triangles--5 around a point}), cube ({\em squares--3
around a point}), and dodecahedron ({\em pentagons--3 around a point}) have a
sum of angles around a vertex of $< 360^\circ$.  A plane can be tiled by
{\em squares--4 around a point} to make a checkerboard pattern, where the
sum of angles around each point is therefore $4 \times 90^\circ$, or
$360^\circ$ degrees---this is a surface of zero curvature.  Johannes
Kepler recognized that the three long-known planar networks,
{\em triangles--6 around a point}, {\em squares--4 around a point}, and
{\em hexagons--3 around a point}, were also regular polyhedra but with an
infinite number of faces. (Kepler in addition allowed regular star
polygons crossing through each other to count, creating regular
stellated polyhedra. With slightly more lenient rules you can find
additional interesting structures.)

In early 1965, found 7 regular spongelike polygon networks having a
sum of angles around a point $> 360^\circ$, with an infinite number of faces
and an infinite number of holes:  {\em triangles--8 around a point},
{\em triangles--10 around a point}, {\em squares--5 around a point}, {\em squares--6
around a point}, {\em pentagons--5 around a point}, {\em hexagons--4 around a
point}, and {\em hexagons--6 around a point}.  I called these
pseudopolyhedrons (Fig.~1), after the pseudosphere which is a
surface of constant negative curvature encountered in the
non-Euclidean geometry of Nikolai Lobachevsky and Janos Bolyai.

\begin{figure}[t!]
\centering
\includegraphics[width=12.6cm]{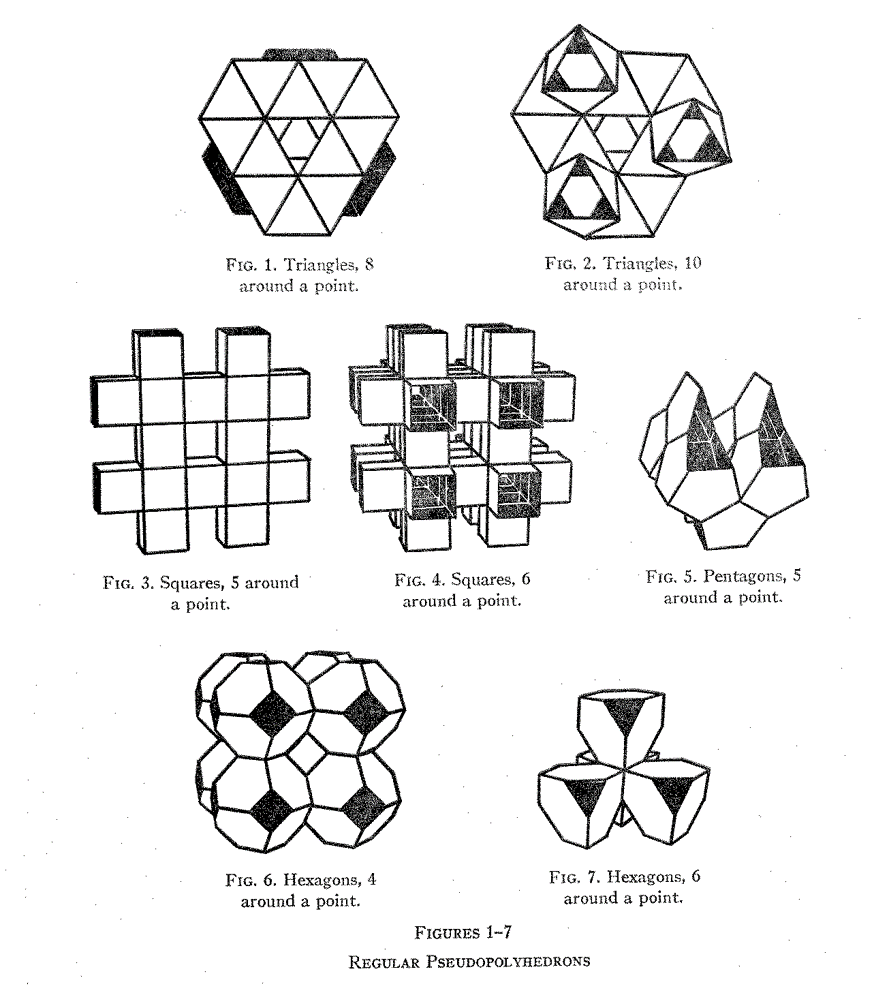}
\caption{Regular Pseudopolyhedrons from Gott (1967).\label{fig:f1}}
\end{figure}

After I got to Harvard, math instructor Tom Banchoff, later famous for
his professional friendship with Salvador Dali, encouraged me to
submit a paper on pseudopolyhedra to the American Mathematical
Monthly, which I did.  The referee’s report was quite positive, but
noted that three of my polygon networks had been discovered before.
The reference was to a paper, which I had never heard of, by
H. S. M. Coxeter (1937).  It described how the first of these
figures---{\em squares--6 around a point}---was found in 1926 by John
Petrie, who also discovered {\em hexagons--4 around a point} (the one
I found first when I was 18). Petrie gets credit for discovering this
entire class of figures.  Coxeter himself discovered {\em hexagons--6
  around a point}.  Petrie and Coxeter did this work in 1926 when both
were 19.  In addition to demanding that the configuration of polygons
around each vertex be identical (as I did), their criteria for
regularity also demanded that the dihedral angles between all adjacent
pairs of faces also be equal and the vertex have rotational symmetry.
With those conditions they were able to prove that the three examples
they found were the only regular figures of this type.  They called
them regular skew polyhedra.  I was happy to add the Coxeter/Petrie
reference.  My paper was still publishable, the referee said, because
I had discovered four new pseudopolyhedrons.  I still required the
configuration of polygons around each vertex to be identical, but
allowed the dihedral angles between adjacent faces to vary.  Some had
some dihedral angles of $180^\circ$, for example. I rediscovered all
three structures discovered by Petrie and Coxeter as well as finding
four new ones allowed by my more lenient rules.  My paper appeared in
print in 1967.  It was my first published scientific paper. My
pseudopolyhedrons were, along with the Coxeter/Petrie regular skew
polyhedra, were included in the {\em Penguin Dictionary of Curious and
  Interesting Geometry} (Wells, 1991).  In my paper I used an anglicized plural
form, pseudopolyhedrons, thinking this would be modern, but the Latin
plural has remained healthy over the past half century, so I will
surrender and here refer to them in the plural as
pseudopolyhedra. These are today sometimes also called infinite
polyhedra, spongelike polyhedra, or infinite skew polyhedra.  These
have an infinite number of faces and therefore belong to the set of
apeirohedra (along with the regular planar networks and cylindrical
networks).

When Siobhan Roberts wrote her definitive biography of Coxeter, {\it
  King of Infinite Space}, in 2006, I was happy to contribute my story
of the astronomical applications these figures later had in
understanding the distribution of galaxies in space. In the early
1980’s there were two schools of thought about how galaxies were
clustered in space. The American school, headed by Jim Peebles,
maintained that there was a hierarchical pattern of clusters of
galaxies floating in a low density void, like isolated meatballs in a
low density soup. The Soviet school, headed by Yakov Zeldovich,
maintained that galaxies formed on a giant honeycomb punctuated by
isolated voids.  I realized neither model was consistent with the new
theory of inflation, which showed that fluctuations in density in the
early universe were produced by random quantum fluctuations. In this
case, the regions of above-average and below-average density should
have the same topology. This could occur with a spongelike topology,
which divided space into two equivalent parts. I knew this because it
occurred in some of my pseudopolyhedra ({\em triangles--10 around a
  point}, {\em pentagons-5 around a point}, {\em squares--6 around a
  point}, {\em hexagons--4 around a point}, {\em hexagons--6 around a
  point}). We showed (Gott, Melott, and Dickinson 1986) that the
spongelike initial conditions required by inflation would grow under
the influence of gravity into a spongelike structure of galaxy
clusters connected by filaments of galaxies, with low density voids
connected by tunnels, a structure now known as the cosmic web and
verified by many surveys.  I tell the story of this discovery in my
book {\it The Cosmic Web} (Gott, 2016).  I gave an invited lecture on
this at the Royal Institution which can be seen on YouTube.

Additional regular pseudopolyhedra have been discovered by
crystallographer A. F. Wells: including {\em triangles--7 around a point};
{\em triangles--9 around a point}; and {\em triangles--12 around a point}.
These are illustrated in Wells's 1969 paper and his 1977 book, {\it Three
Dimensional Nets and Polyhedra}.  Wells, like me, did not demand equal
dihedral angles between adjacent faces.  All are spongelike with an
infinite number of faces and an infinite number of holes. Melinda
Green rediscovered my {\em pentagons--5 around a point}, and has
illustrated many pseudopolyhedra (see references).  Avraham Wachmann,
Michael Burt, and Menachem Kleinman (abbreviated WBK) have discovered
many semi-regular spongelike polyhedra, composed of polygons of more
than one kind, for example, two squares and two hexagons around each
point.  (But they failed to find {\em pentagons--5 around a point}.) These
are to the Petrie/Coxeter/Gott/Wells pseudopolyhedra as the
Archimedean polyhedra are to the 5 classic Platonic polyhedra and are
illustrated in their 1974 book {\em Infinite Polyhedra}. WBK also allow
different dihedral angles between faces as I and Wells did, and in
addition allowed networks to contain pairs of mirror vertices, ones
where the arrangement of polygons was congruent only under mirror
reflection. Envelope polyhedra containing such mirror vertices will be
discussed in the second half of the Appendix.

\section{Finding Envelope Polyhedra---Dihedra} \label{sec:finding}

I thought of envelope polyhedra in 1991, while visiting Aspen,
Colorado to attend a seminar on cosmology, to talk about my two-moving
cosmic string solution in general relativity which allowed time travel
to the past.  I was thinking about the classic regular polyhedra as
approximations to a sphere.  For this reason, there are polyhedral
maps of the Earth.  Perhaps the most famous and successful is the
Gnomonic Cahill Butterfly map.  It maps the Earth onto a regular
octahedron.  Then one unfolds the 8 triangular faces in a butterfly
pattern sitting on a plane.  It shows relatively low distortion but
has a number of ``interruptions.''  Another successful polyhedron map
was invented by the famous architect Buckminster Fuller, the inventor
of the geodesic dome.  Fuller mapped the Earth onto a 20-sided
icosahedron.  He unfolded the 20 triangular faces to sit on the plane.
I remembered that there was also a conformal projection invented by
Emile Guyou in 1887, which maps the two hemispheres of the Earth onto
two squares sitting side by side.  If one folds the two squares
together as one folds a billfold closed, and seals them one would
create an envelope with the Western Hemisphere mapped on the front and
the Eastern Hemisphere on the back.  This envelope has two square
faces taped together along their edges back to back.  I realized this
is also a polyhedron.  It has 2 square faces (the front and back of
the envelope), 4 edges (which form the edges of the envelope) and 4
vertices (which form the four corners of the envelope.  This obeys
Euler's rules for convex polyhedra, namely that the number of faces
minus the number of edges plus the number of vertices equals 2: $F - E +
V = 2$.  This envelope has $F = 2$, $E = 4$, $V = 4$, so $2 - 4 + 4 =
2$.  It also obeys Descartes' rule that the sum of the angle deficits
in a convex polyhedron must be $720^\circ$.  Only 2 faces come together at a
vertex: the front of the envelope and the back.  This is a polyhedron
we would designate as {\em squares--2 around a point}.  This has a Schl\"afli
symbol \{4,2\} and in the WBK nomenclature would be designated
$4^2$. Imagine an ant tethered to a vertex (one of the four corners) with
a tiny string.  It stays at a constant distance from the vertex as it
circles it.  An ant crawling around this vertex would traverse an
angle of $90^\circ$ on the front square, then go over an edge and start
crawling on the back square through another $90^\circ$.  So the total angle
at the vertex is $90^\circ + 90^\circ = 180^\circ$.  This is $180^\circ$ less
than we would get circling a point on a plane which is $360^\circ$,
giving an angle deficit for this vertex of $180^\circ$.  There are 4
such vertices or corners, and so the total angle deficit is $4 \times
180^\circ = 720^\circ$, just as Descartes would have figured.  This is
of course a polyhedron with zero volume, which is why the ancients did
not count it.  In general, one has envelope polyhedra which are
{\em $N$-gons--2 around a point}, for all $N \geq 3$.  Each has 2 Faces:
{\em triangles--2 around a point}, {\em squares--2 around a point}, {\em pentagons--2
  around a point}, {\em hexagons--2 around a point}, and so forth.

Years
later, I found that this was just a rediscovery on my part of dihedra
which have already been accepted as polyhedra for some time. (Coxeter,
1937) mentioned dihedra, for example (but did not include them in his
lists of regular polyhedra).  They all have the topology of a
sphere. (One can therefore make a conformal map of the Earth on two
hexagons: one covering the northern hemisphere, one covering the
southern hemisphere. This answers in the affirmative the gamer’s
perennial question: can the sphere be tessellated with identical
geodesic hexagons? Yes, with two. Each has six geodesic sides, six $60^\circ$ geodesic arcs, along the equator.) There are an infinite number of
envelope polyhedra with 2 faces, $N$ edges, and $N$ vertices.  These all
satisfy the $F - E + V = 2$ rule for convex polyhedra, in the most
transparent way possible. They also satisfy the Coxeter-Petrie
condition that all dihedral angles be equal (in this case $0^\circ$) and that
the vertex figure should have rotational symmetry (in this case
$n = 2$). These are not new. But they are just a subset of the larger class
of envelope polyhedra that are the subject of this paper. From my work
on pseudopolyhedra I already knew in 1991 that if one allowed dihedral
angles of $0^\circ$, there would be many new envelope polyhedra of zero and
negative curvature, both finite and infinite with $360^\circ$ around a vertex
and $>360^\circ$ around a vertex. And this would make for many additional
interesting structures.  I started adding these to my list of envelope
polyhedra. If dihedra can have dihedral angles between their two faces
of $0^\circ$, then this should be allowed for dihedral angles in general.

\section{A Wealth of Envelope Polyhedra} \label{sec:wealth}

I was used to polyhedra approximating surfaces that divided space into
two regions: the inside and the outside in the case of finite
polyhedra, into two regions in the case of the plane tessellations:
above the plane and below the plane, and into two spongelike
interlocking regions in the case of the pseudopolyhedra.  Envelope
polyhedra, do not divide space into two parts.  They have some
dihedral angles (angles between two faces) that are $0^\circ$.  But
they do represent a surface an ant could crawl over.  I don't allow
the surface to cross itself (I am not considering starred
polyhedra). A regular envelope polyhedron has faces that are regular
polygons, and the arrangement of polygons touching each vertex must be
identical.  Some of the angles between adjacent faces in envelope
polyhedra will be $0^\circ$---those faces will be back to back as in
an envelope. This leads to a wealth of forms.

\section{Finite Envelope Polyhedra with Holes} \label{sec:finite}

See Fig.~2, {\it squares--6 around a point}. This and following
photographs are stereo pairs with the left-eye view on the left and
the right-eye view on the right. Place your nose on the page and the
left eye view will be in front of your left eye and the right eye view
will be in front of your right eye. You will see a blurry 3D view,
slowly back away from the page and the fused central image will come
into clear focus, with side images to the left and right.  These may
also be viewed with a standard stereograph viewer. The models in
Figs.~2--5 are made from Polydron plastic polygons (from
polydron.com), ignore their slightly serrated edges which allow them
to be hinged together to make polyhedra.

\begin{figure}[t!]
\centering
\includegraphics[width=12.6cm]{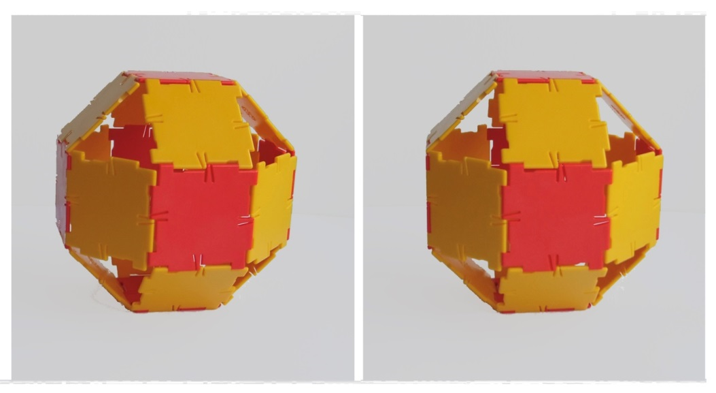}
\caption{The finite envelope polyhedron, {\em squares--6 around a point}.\label{fig:f2}} 
\end{figure}

{\em Squares--6 around a point}, as illustrated in Fig.~2 starts with a
rhombicuboctahedron with 18 square faces and 8 triangular faces.
Remove the triangular faces to create triangular holes into the
interior.  This leaves an envelope polyhedron that is hollow with 18
exterior square faces and 18 interior square faces. The polygonal
faces thus partially envelop a hollow void.  Thus, the name envelope
polyhedron seems appropriate here too.  Each vertex connects 3
exterior and 3 interior squares.  An ant crawling around a vertex
would visit all 6 of these squares.  [The pseudorhombicuboctahedron, a
  Johnson polyhedron, where the top in Fig.~2 is twisted relative to
  the bottom by $45^\circ$ provides an alternate way to start the
  construction and gives an alternate form (2)].  All vertices are
congruent with a sum of angles of $6 \times 90^\circ > 360^\circ$
around each, so this approximates a negatively curved surface. It has
36 faces (18 exterior plus 18 interior faces), $6 \times 4 = 24$
interior edges, $6 \times 4 = 24$ exterior edges, and $8 \times 3 =
24$ edges on the 8 triangular holes connecting the interior and
exterior for a total of 72 edges, as well as $8 \times 3 = 24$
vertices.

Note that just as there are pairs of faces (exterior and
interior) that are back to back, there are also pairs of edges that
are back to back. But there are no vertices back to back. Vertices
connect the exterior and interior faces.  Cut one triangular hole in
the rhombicuboctahedron and it looks like a bowl with a small mouth
which can be distorted into a disk and fattened into a sphere. Each of
the additional 7 triangular holes you cut creates a doughnut hole and
creates a handle on the sphere. So the genus of this figure (number of
doughnut holes) is 7. It obeys the rule $F - E + V = 2 \times (1 -
\mbox{genus}) = 36 - 72 + 24 = -12$. This gives it a genus of 7 and a
topology equal to that of a sphere with 7 handles. It has negative
curvature and so could be considered a finite envelope
pseudopolyhedron.  Note that the faces around each vertex must create a
single surface around each vertex that an ant tethered to the vertex
could visit. Thus two cubes touching at a point would not be
considered {\em squares--6 around a point} because an ant tethered to the
vertex would circle the vertex visiting 3 squares on one cube and
complete her circuit and return to where she started visiting only
those 3 squares and never visit the other cube. Likewise, a cube with
a square fin attached at an edge would not count as part of a
{\em squares--5 around a point} structure because although one ant tethered
at the vertex would visit 5 outside square faces, another ant tethered
on the inside of the cube would visit only 3 squares, again creating
multiple surfaces with different angle deficits at a single vertex. We
are not allowing such multiple surfaces. The envelope polyhedra is one
continuous two-dimensional surface.

Here is another interesting
example: {\em octagons--4 around a point}.  Get an empty cubical cardboard
box.  Cut off each of its corners with a saw.  This will make a
truncated cube.  Each of the 6 square faces of the cube will have its
corners cut off, becoming an octagon.  Cutting the corners of the
cubical box will create 8 triangular holes where the 8 corners of the
box used to be.  This has 12 octagonal faces (6 inside, 6 outside), 48
edges (12 outside edges, 12 inside edges, and 24 edges on the 8
triangular holes connecting the inside and outside), and 24 vertices
(3 on each of the triangular holes). $F - E + V = 12 - 48 + 24 = -12$.
That’s $2 \times (1 - \mbox{genus})$ as predicted.  This has a genus
of 7.  It has the same number of triangular openings as the {\em squares--6
around a point} envelope polyhedron in Fig.~2 so it also has a genus
of 7.  Some dihedral angles between adjacent octagons are $90^\circ$ (when
both are either outside or both inside) and some are $0^\circ$ (when one is
inside and one is outside).

An ant tethered to one of the vertices
will circle it by crawling over two outside-facing octagons and two
inside-facing octagons as it goes around the outside and then through
the triangular hole to visit the inside, making 4 octagons around each
vertex.  The interior angle in an octagon is $135^\circ$.  So the
total angle the ant traverses circling the vertex is $4 \times
135^\circ = 540^\circ > 360^\circ$. This is a negatively curved
surface with the curvature all concentrated at the vertices.  It even
looks like a saddle-shaped surface.  One can imagine an ant sitting on
one of the vertices and draping its little hind legs over each side,
one inside and one outside, like he was riding a horse.  This
negatively curved finite polyhedron, like {\em squares--6 around a point} in
Fig.~2, is multiply connected.  But all of my original multiply
connected pseudopolyhedra were infinite.

{\em Hexagons--4 around a
point}. Take an octahedron shaped box and saw off its corners. You will
be left with a truncated octahedron shape made of hexagons with 6
square holes. See Fig.~3. An ant circling a vertex will traverse two
hexagons on the exterior, before visiting another two hexagons on the
interior to make {\em hexagons--4 around a point}.

\begin{figure}[t!]
\centering
\includegraphics[width=12.6cm]{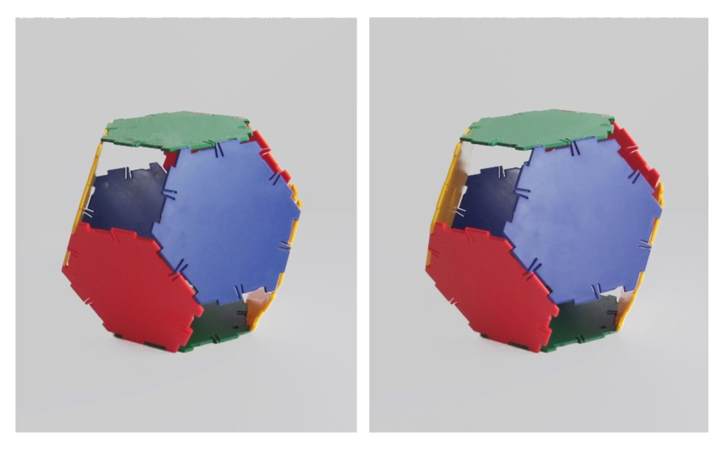}
\caption{Envelope polyhedron, {\em hexagons--4 around a point}.\label{fig:f3}} 
\end{figure}

\begin{figure}[t!]
\centering
\includegraphics[width=12.6cm]{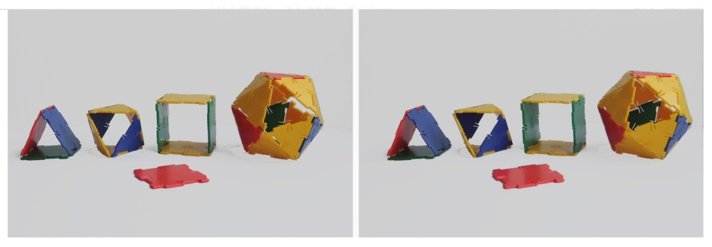}
\caption{Finite envelope polyhedra: Front, {\em squares--2 around a point}; Back row, left to right, {\em squares--4 around a point}, {\em triangles--6 around a point}, {\em squares--4 around a point}, {\em triangles-8 around a point}.\label{fig:f4}} 
\end{figure}

Take a triangular prism and remove the top and bottom triangles, to
create the envelope polyhedron {\em squares--4 around a point}. Take a
regular octahedron and remove two opposite faces, to create the
envelope polyhedron {\em triangles--6 around a point}. Take a cube and
remove two opposite faces to create the envelope polyhedron {\em squares--4
around a point}. These have the topology of a doughnut, and $360^\circ$  and
zero curvature at the vertices. Take an icosahedron box and remove
four triangles to create the envelope polyhedron {\em triangles--8 around a
point}. These are illustrated in Fig.~4.

\section{Infinite Envelope Polyhedra} \label{sec:infinite}

{\it Squares--10 around a point}. The hollow rhombicuboctahedron envelope
polyhedron shown in Fig.~2 fits perfectly in an imaginary cube.
Stack cubes like this to fill space with the rhombicuboctahedra glued
back to back at a single square (red in Fig.~2) in a repeating
pattern.  Each vertex connects two of these rhombicuboctahedra, and as
an ant circles the vertex, it travels on an exterior square face of
the first, then three interior faces of the first, then another
exterior face of the first, an exterior face of the second, three
interior faces of the second, and another exterior face of the
second---giving 10 squares around the vertex.

{\it Octagons--8 around a point}.  Take an infinite number of the
truncated cubical cardboard boxes ({\em octagons--4 around a point})
discussed above and fill all of space with them by stacking them like
cubical boxes in a warehouse.  We will be gluing the boxes together.
Two adjacent boxes will be glued together on their outside octagons.
These will disappear from the surface.  All that will be left are the
inside octagons of each cube.  Four boxes (numbered 1, 2, 3, 4) will
fit together at their edges.  A vertex will connect 8 octagons.  The
vertex is at the end of an edge that connects two inside octagons of
cube 1.  Cubes 2, 3, and 4 also come together at this vertex.  An ant
tethered to the vertex will circle the vertex by first traversing the
two inside octagons of cube 1, then it will cross an edge of a
triangular hole, to enter cube 2 and traverse two of its inside
octagons, then enter cube 3 traversing two of its interior octagons,
before entering cube 4 and traversing two of its interior octagons as
she returns to the place she started.  This makes a complicated saddle
shaped surface that goes up and down, up and down, up and down, and up
and down.  It is a saddle with a four-fold symmetry, like a +.  It is
a saddle a horse could sit on, hanging its legs down, one each
into boxes 1, 2, 3, 4.  It is cradled by the four triangular holes.
The angle around each vertex is $8 \times 135^\circ = 1080^\circ$.
This is an infinite polyhedron: it has an infinite number of faces, an
infinite number of cubical cells with triangular holes.  A fly could
fly through the whole structure, visiting any cell he wanted.

\begin{figure}[t!]
\centering
\includegraphics[width=12.6cm]{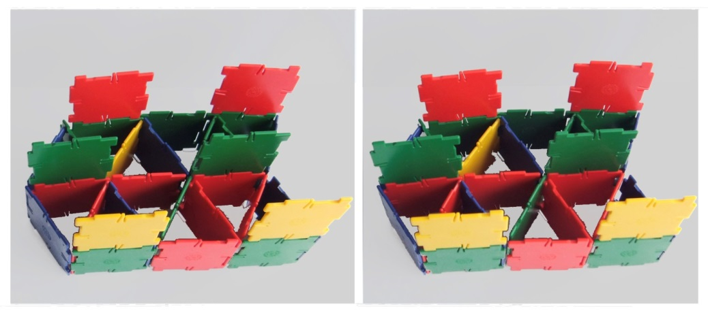}
\caption{{\em Squares--14 around a point}.\label{fig:f5}} 
\end{figure}

(Note:
Pellicer and Schulte (2009) and Schulte and Weiss (2016) considered
skeletal polyhedra. These had skeletal faces consisting only of edges,
based on the idea of Gr\"unbaum (1994) that one could consider polyhedra
consisting of polygons with edges only without membranes spanning them
to form a surface. These form additional polyhedra, having skew
polygons, and even infinite helical polygons. In these polyhedra
skeletal polygons come together two at an edge. But Pellicer and
Shulte (2009) also consider regular polygonal complexes where more
than two skeletal faces meet at an edge. If we removed all the faces
from {\em octagons--8 around a point}, we would have 4 skeletal octagons
meeting at an edge. But to be clear, we are not considering such
skeletal structures here---envelope polyhedra all have surfaces. In
{{\em octagons--8 around a point} we have 8 octagons and 8 edges meeting at a
point. In the surface, two octagons always share an edge. Edges can be
placed back to back, of course, in the structure, just as polygons can
be placed back to back. These constitute, in the surface, separate
polygons and separate edges which an ant would encounter circling the
vertex.)

{\em Triangles--18 around a point}.  (1) Start with the
pseudopolyhedron triangles--10 around a point shown in Fig.~1, then
remove triangles in the planes of triangles that are connected to
other triangles in the same plane (with dihedral angles of $180^\circ$) along
each of their three edges.  This creates triangular holes in the
planes of triangles.  Each vertex used to have 10 triangles around it
on one convoluted surface; take one triangle away and that leaves 9
triangles, but the hole allows the ant to traverse to the opposite
side of the plane and visit the 9 opposite sides of these triangles to
make a total of 18 around a point.  The angle around each vertex is
$18 \times 60^\circ = 1080^\circ$. There are left- and right-handed versions of this.

{\em Squares--12 around a point}.  Make a single layer of triangular prisms
that cover a plane.  Delete the triangular tops and bottoms of the
prisms, leaving a triangular grid of cells.  Six triangular prisms
meet at a vertex, and each has two interior squares that touch the
vertex, giving 12 squares around a point.

{\em Squares--14 around a point}.
Take the previous structure and add east-west fins of back-to-back
squares above this single layer of cells in a dashed pattern.  See
Fig.~5. This way each vertex has the 12 squares around a point that
it had before plus the two squares in the fin.  These fins are then
connected to another plane of {\em squares--12 around a point} triangular
cell layer above it.  Layers of fins and cells alternate vertically
forever.  Each of the 14 squares has a vertex angle of $90^\circ$, so the
total angle around a vertex is $(14 \times 90^\circ) = 1260^\circ$, a
trick for an Olympic snowboarder.  That is $3\nicefrac{1}{2}$ rotations.  That
gives it an angle excess of $900^\circ$ (above $360^\circ$), the most found
for any polyhedron so far.

There are many other envelope polyhedra.  The Appendix gives a list of
those I have found so far and a summary of how these additional ones
are constructed.  In general they are found by: (a) deleting faces
from regular polyhedra to create polygonal holes (as illustrated by
some of the examples in Fig.~4); (b) deleting polygons from
semi-regular polyhedra containing faces of different types (like the
triangular prism in Fig.~4) and the Archimedean polyhedra to create
polygonal holes; (c) deleting polygons from regular and semi-regular
planar networks and cylindrical networks; (d) deleting polygons from
regular pseudopolyhedra like those shown in Fig.~1, and from
semi-regular pseudopolyhedra like those shown in WBK, and (e) Creating
planar networks of open-ended prisms with fins (like Fig.~5) and
without fins, based on regular and semi-regular (Archimedean) planar
networks.

\section{Summary} \label{sec:summary}

Polyhedra fall into three groups depending on whether the sum of face
angles around a vertex are $< 360^\circ$, $= 360^\circ$, or $>
360^\circ$ (i.e., whether the curvature at each vertex is positive,
zero, or negative. Dihedral angles do not all have to be
identical. But the arrangement of polygons around each vertex must be
identical.  Envelope Polyhedra, which we introduce here, just have
some dihedral angles which are $0^\circ$.  Using the nomenclature of WBK,
where {\em $N$-gons--$M$ around a point} with Schl\"afli symbol $\{N, M\}$ are
designated $N^M$, we find the following structures so far:

\begin{tabular}{ll}
\\
  $<360^\circ$\\
  &$3^3$, $3^4$, $3^5$, $4^3$,$5^3$ (The classical Platonic Polyhedra)\\
&$N^2$ where $N \geq 3$ (Finite Envelope Polyhedra)\\
\end{tabular}

\begin{tabular}{ll}
  $=360^\circ$\\
  &$3^6$, $4^4$, $6^3$ (Infinite Planar and/or Cylindrical Tessellations)\\
  &$3^6$, $4^4$ (Finite Envelope Polyhedra with Toroidal Geometries)\\
  &$3^6$, $4^4$ (Infinite Envelope Polyhedra with Filmstrip Geometries)\\
\end{tabular}

\begin{tabular}{ll}
$>360^\circ$\\
  &$3^7$, $3^8$, $3^9$, $3^{10}$, $3^{12}$, $4^5$, $4^6$, $5^5$, $6^4$, $6^6$
  (Infinite Pseudopolyhedra/Skew Polyhedra)\\
  &$3^8$, $3^{10}$, $4^6$, $6^4$, $8^4$, $10^4$
  (Finite Envelope Pseudopolyhedra)\\
  &$3^{10}$, $3^{12}$, $3^{14}$, $3^{18}$, $4^{6}$, $4^{10}$, $4^{12}$, $4^{14}$,
  $6^8$, $8^4$, $8^8$
  (Infinite Envelope Pseudopolyhedra)\\
  \\
\end{tabular}

These are summarized in Table 1.

\begin{table}[h!]
  \centering
  \begin{tabular}{|l|l|l|l|l|l|l|l|l|l|}
    \hline
    $N^2$ $\diamond$ & $N^3$ & $N^4$ & $N^5$ & $N^6$ & $N^7$ & $N^8$ & $N^9$ & $N^{10}$ & $N^{11}$ \\
    \hline
    $12^2$ $\diamond$ & $12^3$ & $12^4$ $\circ$ & $12^5$ & $12^6$ & $12^7$ & $12^8$ & $12^9$ & $12^{10}$ & $12^{11}$ \\
    $11^2$ $\diamond$ & $11^3$ & $11^4$ & $11^5$ & $11^6$ & $11^7$ & $11^8$ & $11^9$ & $11^{10}$ & $11^{11}$ \\
    $10^2$ $\diamond$ & $10^3$ & $10^4$ $\diamond$& $10^5$ & $10^6$ & $10^7$ & $10^8$ & $10^9$ & $10^{10}$ & $10^{11}$ \\
    $9^2$ $\diamond$ & $9^3$ & $9^4$ & $9^5$ & $9^6$ & $9^7$ & $9^8$ & $9^9$ & $9^{10}$ & $9^{11}$ \\
    $8^2$ $\diamond$ & $8^3$ & $8^4$ $\diamond\circ$ & $8^5$ & $8^6$ & $8^7$ & $8^8$ & $8^9$ & $8^{10}$ & $8^{11}$ \\
    $7^2$ $\diamond$ & $7^3$ & $7^4$ & $7^5$ & $7^6$ & $7^7$ & $7^8$ & $7^9$ & $7^{10}$ & $7^{11}$ \\
    $6^2$ $\diamond$ & $6^3$ $\bullet$ & $6^4$ $\diamond\bullet$& $6^5$ & $6^6$ $\bullet$ & $6^7$ & $6^8\circ$ & $6^9$ & $6^{10}$ & $6^{11}$ \\
    $5^2$ $\diamond$ & $5^3$ $+$ & $5^4$ & $5^5$ $\bullet$ & $5^6$ & $5^7$ & $5^8$ & $5^9$ & $5^{10}$ & $5^{11}$ \\
    $4^2$ $\diamond$ & $4^3$ $+$ & $4^4$ $\diamond\bullet\circ$ & $4^5$ $\bullet$ & $4^6$ $\diamond\bullet\circ$& $4^7$ & $4^8$ & $4^9$ & $4^{10}$ $\circ$ & $4^{11}$ \\
    $3^2$ $\diamond$ & $3^3$ $+$ & $3^4$ $+$ & $3^5$ $+$ & $3^6$ $\diamond\bullet\circ$& $3^7$ $\bullet$ & $3^8$ $\diamond\bullet$& $3^9$ $\bullet$ & $3^{10}$ $\bullet\circ$ & $3^{11}$ \\
    \hline

  \end{tabular}

  \begin{tabular}{|l|l|l|l|l|l|l|}
    \hline
    $4^{12}$ $\circ$ & $4^{13}$ & $4^{14}$ $\circ$& $4^{15}$ & $4^{16}$ & $4^{17}$ & $4^{18}$\\
    $3^{12}$ $\bullet\circ$& $3^{13}$ & $3^{14}$ $\circ$& $3^{15}$ & $3^{16}$ & $3^{17}$ & $3^{18}$ $\circ$\\
    \hline

  \end{tabular}

\begin{tabular}{rcl}
  $+$ & = & Finite  Polyhedron\\
  $\diamond$ & = & Finite Envelope Polyhedron\\
  $\bullet$ & = & Infinite Planar \& Cylindrical Networks, or Pseudopolyhedra\\
  $\circ$ & = & Infinite Envelope Polyhedron\\
\end{tabular}
  
  \caption{Nomenclature follows that of WBK: For example {\em squares--6
    around a point} appears as $4^6$. This has a Schl\"afli symbol
    \{4,6\}. The symbols $\diamond$ and $\circ$ in that box indicate that both
    finite and infinite envelope polyhedra of this type exist; $\bullet$
    indicates that an infinite pseudopolyhedron of this type exists. $N$
    represents all integers $N \ge 3$.}
  
  \label{tab:table1}
\end{table}

A number of these have several
geometrical forms as described in the Appendix.  Open symbols refer to
envelope polyhedra which have some dihedral angles equal to $0^\circ$: the
symbol $\diamond$ indicates a finite number of sides and refers to polyhedra
with a finite number of faces, $\circ$ is a circle and refers to envelope
polyhedra with an infinite number of faces. Closed symbols refer to
Polyhedra: $+$ (which are finite), and $\bullet$ which are either infinite
Planar or Cylindrical Networks, or infinite Pseudopolyhedra which have
no dihedral angles equal to $0^\circ$.

\section{Envelope Polyhedra with Mirror Vertices} \label{sec:mirror}

\begin{figure}[b]
\centering
\includegraphics[width=12.6cm]{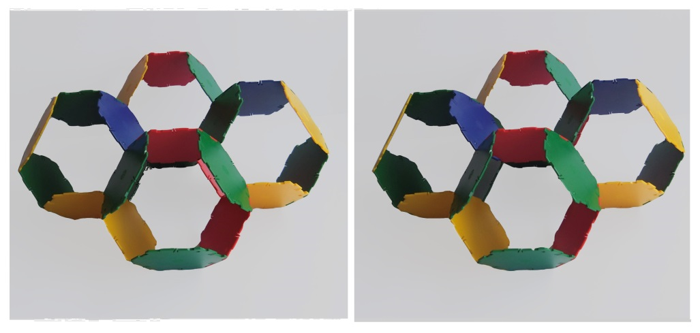}
\caption{{\em Hexagons--6 around a point} with mirror image
  vertices.}
\label{fig:f6}
\end{figure}

Below, and in the Appendix, we also list additional infinite
structures ($3^{12}$, $4^6$, $4^8$, $4^{10}$, $6^6$) with mirror
vertices. All the previous polyhedra we have discussed have vertices
that are identical. But one may obtain a considerable number of
additional structures by allowing vertices and their mirror images to
be considered as congruent. One example starts with the
pseudopolyhedron {\em hexagons--4 around a point} shown in Fig. 1. Now
remove the same set of two opposite faces on each truncated octahedron
box. You get a repeating pattern of {\em hexagons--6 around a point}:
$6^6$, but with pairs of mirror vertices whose vertex figures are not
identical (except under mirror reflection) which is shown in Fig.~6.
Because it includes the mirror image vertices, this structure $6^6$ is
not included in Table~1 which includes only structures with vertices
that are identical.

Also starting with {\em squares--6 around a point} in Fig.~1, on can delete
horizontal squares that are tops and bottoms of cubes open in the
east-west direction to get {\em squares--10 around a point}, which also
includes mirror vertices (with mirror symmetric vertex figures).

Many additional such structures with mirror vertex figures can be
obtained by deleting polygons from structures in WBK. I should mention
that WBK consider mirror image vertices as congruent and many of their
semi-regular pseudopolyhedra have vertices that are congruent only
under mirror reflection. In plane geometry we are used to considering
a triangle and its mirror image as congruent, but in three dimensions
they are actually identical since they can be superimposed by rotation
in the third dimension. Chemists know that in 3D, glucose comes in
distinct mirror image right-handed (\textsf{D}-glucose = dextrose) and
left-handed (\textsf{L}-glucose) forms. Dextrose can be digested by
the body to produce energy while \textsf{L}-glucose cannot. In 3D
solid Euclidean geometry mirror image vertices cannot be
superimposed. So some might regard those with mirror image vertices as
less regular. Therefore I have separated the envelope polyhedra into
two classes, those whose vertices have identical vertex figures,
listed in the previous sections and the first part of the Appendix,
and those containing mirror image vertex figures listed in this
section and in a separate section of the Appendix.

No doubt there are
more envelope polyhedra yet to be discovered.  I leave that as a
challenge for interested readers---to add to the list. Surprisingly,
pseudopolyhedra eventually had an application in astrophysics. Will
any of these envelope polyhedra eventually find applications, in
chemistry, biology, art, or architecture perhaps?  Time will tell.

\acknowledgments

I thank Robert J. Vanderbei for 3D photography of the models, and
Wesley N. Colley for help in preparation of the manuscript.

\appendix
\section{Additional Envelope Polyhedra}

\begin{itemize}
\item
{\em Triangles--6 around a point.}  (1) Start with an $N$-gon antiprism ($N
\ge 3$) and remove the top and bottom $N$-gons.  This leaves a band of
triangles connecting them, like the corrugated circumference of a drum
without the top and bottom.  Each vertex links 3 outside triangles
with 3 interior triangles.  These have $2N$ faces with the topology of a
doughnut.  The $N = 5$ case is an icosahedron with its top and bottom caps
of 5 triangles each deleted, leaving a ring of 10 exterior triangles
and 10 interior triangles.  (2) A straight filmstrip of triangles that
is infinitely long.  They alternate between triangles with points up
and points down.  Each vertex links 3 triangles on one side with 3
triangles on the other side.  These vertices have zero curvature and
$360^\circ$ around a vertex. (3) Such an infinite filmstrip may be twisted to
form a helix, there are left- and right-handed versions of these.

\item
{\em Triangles--8 around a point.}  (1) Start with a hollow snub cube, one
of the Archimedean polyhedra where 4 triangles and one square meet at
each vertex, and delete the square faces.  This leaves 32 exterior
triangular faces, 32 interior triangular faces, and 6 square holes
connecting the interior and exterior.  There are left and right handed
versions of this.  (2) Start with a hollow snub dodecahedron, another
Archimedean polyhedron, and delete the pentagonal faces.  This leaves
80 exterior triangular faces, 80 interior triangular faces, and 12
pentagonal holes connecting the interior and exterior.  This has the
topology of a sphere with 11 handles.  There are both left and right
handed versions of this.

\item
{\em Triangles--12 around a point.}  (1) Start with {\em squares--8
  around a point}, a single plane layer of cubical boxes with open
ends. This is a checkerboard of cubical boxes with open tops and
bottoms. The vertical sides of the boxes are oriented either east-west
or north-south. All the boxes have vertical edges. Distort this by
leaning all the vertical edges in the northwest direction until they
are at $45^\circ$ to the vertical, like the diagonal of a vertical
octahedron. The boxes now have rhombus sides that can be tiled by two
equilateral triangles. Each cube is now leaning, with 4 rhombus
interior sides made up of two equilateral triangles each. The rhombus
sides of each open-ended box have angles of $120^\circ$ (where two
triangles meet) and $60^\circ$ where there is one triangle. An ant
circling a vertex on top of the layer would cover 2 interior triangles
in the southeast box, 3 interior triangles in the northeast box, 4
interior triangles in the northwest box, and 3 interior triangles in
the southwest box before returning to where it started, making 12
triangles around a point. Flip the layer of boxes over to see that the
vertices on the bottom of the layer are identical to those on the
top. The total angle around the vertex is $720^\circ$, just as in the
{\em squares--8 around a point} structure we started with. Both are
topologically equivalent. Both structures are flexible.

\item
{\em Triangles--14 around a point}.  (1) Start with the Gott (1967)
pseudopolyhedron {\em triangles--10 around a point}, and pull out for
consideration two adjacent planes of triangles and the tunnels that
connect them.  Each vertex in the plane of triangles will have a
triangular hole in it where the column would have been to connect it
to the next plane of triangles in the {\em triangles--10 around a point}
configuration.  The two planes of triangles, with their triangular
holes are connected to each other by tunnels that are antiprisms
(octahedrons with two missing [top and bottom] opposite faces).  For
example an ant circling each vertex on the top plane would traverse 1
triangle on the bottom of that plane of triangles, 3 triangles that
are part of the exterior of a tunnel, 3 more triangles on the bottom
of the top plane, then through the open triangular hole onto the top
of the top plane, where it would traverse 3 triangles, on the top of
the plane, 3 triangles in the interior of the tunnel, 1 more triangle
on the top of the plane, and it is back where it started. That gives
14 triangles around a point.  There are left and right handed versions
of this. (2) Start with Wells (1977) pseudopolyhedron $3^8$ made of snub
cubes connected at their missing square faces in a cubic array. Left
and right handed snub cubes alternate like black and white cubes in a
3D chessboard. Left handed ones join to right handed ones at the
vertices of the missing squares.  Now
delete from the array of triangles those triangles on left-handed snub
cubes that touch three different missing squares. That leaves
left-handed snub cubes with 8 triangular holes. Each vertex lies on
the corner of a square tunnel connecting a left and right handed snub
cube. It is missing one triangle out of 8 giving its ant a triangular
hole to allow it to visit both sides of the 7 remaining triangles to
give a surface with 14 triangles around a point. There are left and
right-handed versions of this structure (depending on whether one
takes the triangles out of the left-handed, or right-handed snub
cubes. But within each case the vertices are all identical---either all
left-handed or all right-handed.)

\item
{\em Squares--4 around a point}.  (1) Start with an $N$-gon prism (for $N > 2$)
with an $N$-gon top, $N$ square sides, and an $N$-gon bottom, then remove
the $N$-gon top and bottom.  This leaves an envelope polyhedron with 2N
faces, N outside square faces and N inside square faces.  There are an
infinite number of these which all have the topology of a doughnut-like
the sides of a drum with the top and bottom missing.  Each vertex
joins two outside squares and two inside squares.  Thus, the
polyhedron has an angle deficit of $0^\circ$.  For $N = 4$, this is a cube with
the top and bottom missing.  (2) There is also an arrangement where
one has two strips of squares back to back extending in a straight
line like a filmstrip.  This has the topology of an infinite cylinder.
It has zero curvature.  There are also arrangements where the squares
line up back to back in a zig-zag patterns, like a filmstrip with folds
in it back and forth.

\item
{\em Squares--6 around a point.} (1) Start with just the top plane of
Gott's (1967) pseudopolyhedron {\em squares--5 around a point}
($4^5$). This is a single checkerboard plane punctuated by holes so
that three squares on the top surface of the plane surround each
vertex, while three more squares surround it on the bottom surface of
the plane. As the ant circles the vertex it visits 6 squares giving
$4^6$---or the holes can be in a staggered arrangement as in the front
plane of squares in $4^5$ in WBK at the bottom of page 2, or a plane
of squares in $4^5$ in WBK at top right of page 16. (2) A single honeycomb
layer of hexagonal prisms tiling a planar layer, with the tops and
bottoms removed. An ant circles a vertex visiting the two interior
squares of each of three hexagonal cells making {\em squares---6
  around a point}.  This is like Fig. 5 without the fins and with
hexagonal cells instead. (3) Squares and octagons ($4\times8^2$) tile
the plane, to use the nomenclature of WBK. Convert them into a layer
of octagonal and square prisms (cubes) without their tops and
bottoms. Each vertex is surrounded by two interior squares of the cube
and two interior squares of each of two octagonal prisms: 6 squares
around a point.  (4) Triangles and dodecagons ($3\times12^2$) tile the
plane. Similarly, make them into triangular and dodecagonal prisms
without their tops and bottoms. This makes a single layer with 6
squares around a point.

\item  
{\em Squares--8 around a point.} (1) Start with the semi-regular
pseudopolyhedron shown by Wachman, Burt, and Kleinmann (WBK) on page 9
at the top. It looks like triangular prisms joined by cubes to form a
surface with 4 squares and one triangle around each vertex (designated
$4^4\times3$). Simply remove the triangles, to make an infinite
envelope polyhedron with {\em squares--8 around a point}.  With
triangles gone there are 4 squares around each point, with two sides
each which the ant traverses circling the vertex. (2) On the same page
at the bottom are hexagonal prisms joined by cubes in a similar manner
($4^4\times6$), remove the hexagons to get {\em squares--8 around a
  point}. (3) On page 15, is a multilayer structure ($4^5$). Instead
of having two planes of squares connected by tunnels as in (Gott
1967), this has planes of squares connected by tunnels going down and
columns going up to adjoining punctuated planes of squares. Remove
selected squares on the planes of squares (those touching at their
corners 2 columns and 2 tunnels) to create an envelope polyhedron {\em
  squares--8 around a point}. (4) A similar $4^5$ structure on page 16
has a different arrangement of columns and tunnels. Remove selected
squares (on the planes of squares, those touching at their corners 2
columns and 1 tunnel on one plane and 1 column and 2 tunnels on the
next and repeat) to create an envelope polyhedron {\em squares--8
  around a point}.  (5) Make a single layer of cubes that cover a
plane.  From the top they look like a checkerboard.  Delete the tops
and bottoms of the cubes, to leave a checkerboard shaped grid of
cells.  Four cubes join at a vertex, so an ant circling the vertex,
traverses two interior squares of cube 1, then two interior squares of
cube 2, then two interior squares of cube 3, and finally two interior
squares of cube 4 giving 8 squares around a point. (6) Take the
structure $4^4\times6$ on page 40 of WBK. Delete the
hexagons. Equivalently, hexagons, triangles and squares tile the
plane, remove the hexagons and triangles and connect such planes with
triangular prisms (tunnels and columns) with their tops and bottoms
off. This gives {\em squares--8 around a point}. (7) There is a
biplane version of the previous structure where the two planes are
connected by triangular tunnels. (8) Take the filmstrip, {\em
  triangles--6 around a point} (2, listed above) and build triangular
prisms on top of it.  Now delete the triangular faces.  This leaves an
infinite row of interlaced triangular prisms without their tops and
bottoms, giving {\em squares--8 around a point}. (9) Take the
structure $4^5$ on page 50 of WBK, Delete the squares in the planes of
squares whose four corners touch two columns and two tunnels. (10)
Take the structure $3\times4^4$ on page 71 of WBK and remove the
triangles to make {\em squares--8 around a point}. This structure has
hollow shells like that shown in Fig.~2, missing their red squares
joined by cubes missing two opposite squares attaching them where
their two missing red squares used to be. (11) Start with the
structure $3\times4^4$ on page 78 of WBK, remove the triangles to
produce envelope polyhedron $4^8$. An ant at each vertex navigates
both sides of 4 squares as it circles the vertex. This has a “diamond”
structure with a hexagonal ring of squares circling each “carbon” bond
as an axis. (12) Start with the Wells (1977) structure $4^5$ on page
87 of WBK and delete squares which are at the center of a flat cross
shaped pattern. These squares have dihedral angles on all their four
edges equal to $180^\circ$. (13) Start with the Wells (1977) structure
$4^5$ on page 88 of WBK and delete the squares that have 4 dihedral
angles of $215^\circ16^\prime$ at their edges, as indicated in WBK
(square COB in their diagram). These deleted squares sit on the
surfaces of truncated octahedrons, the squares that remain lie on
surfaces of hexagonal prisms. (14) Squares and octagons ($4\times8^2$)
tile the plane, to use the nomenclature of WBK. Make them into
octagonal prisms and square prisms (cubes) without their tops and
bottoms. Add fins (as in Fig.~5) to the tops of the octagonal prism
sides where two octagonal prisms meet. Then repeat vertically to make
alternate layers of prisms and fins. Gives 8 squares around a
point. (15) Triangles and dodecagons ($3\times12^2$) also tile the
plane. Similarly, make them into a single planar layer of triangular
and dodecagonal prisms without their tops and bottoms. Add fins (as in
Fig.~5) to each dodecagonal prism side that attaches to another
dodecagonal prism. This likewise gives 8 squares around a point. (16)
Triangles, squares, and hexagons tile the plane
($3\times4\times6\times4$), turn these into prisms without their tops
and bottoms to make a single planar layer. An ant circling a vertex
would visit 2 interior squares of each of 4 prisms, giving 8 squares
around a point. (17) Triangles and hexagons tile the plane
$(3\times6)^2$. Turn them into prisms without their tops and bottoms
to make a single planar layer. An ant circling a vertex would visit 2
interior squares of each of 4 prisms, giving 8 squares around a point.

\item
{\em Squares--10 around a point}. (1) Start with the structure $3\times4^5$ on page
99 of WBK and delete the triangles. This takes envelope polyhedra
{\em squares--6 around a point} (as shown in Fig.~2) and replaces the
yellow squares with cubes without ends linking the envelope polyhedron
to similar copies of itself. The envelope polyhedron {\em squares--6 around
a point} shown in Fig.~2 fits in an imaginary cube, the added cubes
without ends link it along the imaginary cube’s 12 diagonals to
similar envelope polyhedra in nearby cubes. Circling a vertex the ant
will traverse an outside red square of Fig.~2, two outside squares
of a cube without ends, it will then go through the triangular hole,
visit the two inside squares of that same cube without ends, then the
inside red square, the inside two squares of another cube without
ends, before coming out of the triangular tunnel and traversing the
corresponding two outside squares of that cube without ends before
returning to where it started: 10 squares around the vertex. (2) Start
with a single layer of cubes tiling a plane, remove the tops and
bottoms of the cubes, then add fins to every other east-west square
(similar to Fig.~5). Then repeat vertically to create alternating
layers of planes of open cubes and fins. An ant crawling around each
vertex will visit two interior squares of each of 4 cubes joining each
other, plus the two sides of one fin, giving {\em squares--10 around a
point}. (3) Triangles and hexagons tile the plane ($3^4\times6$). Turn these
into prisms without their tops and bottoms to make a single
layer. Each vertex is surrounded by 5 prisms, each with two interior
faces the ant must visit making it {\em squares--10 around a point}. (4)
Triangles and squares tile the plane ($3^3\times4^2$). Turn these into prisms
without tops and bottoms to make a single layer. Each vertex is
surrounded by 5 prisms, each with two interior faces the ant must
visit so it is also {\em squares--10 around a point}. (5) Triangles and
squares tile the plane ($3^2\times4\times3\times4$). Turn these into
prisms without their tops and bottoms. Similarly, this makes {\em
squares--10 around a point}.

\item  
{\em Squares--12 around a point.} (1) Triangles and squares tile the
plane ($3^3\times 4^2$). Turn these into prisms without tops and bottoms to make
a single layer. Add fins to squares joining two open cubes. Then
repeat vertically to make alternating layers of prisms and fins. An
ant will traverse two sides of a fin, then interior pairs of faces of
5 prisms as it circles each vertex, giving {\em squares--12 around a point}.
(2) Triangles and squares tile the plane
($3^2\times4\times3\times4$). Turn these into prisms without their tops
and bottoms. Add fins above the squares connecting two triangular
prisms. Again the ant traverses two sides of the square fin, then
interior pairs of faces of 5 prisms as it circles each vertex:
{\em squares--12 around a point}.

\item
{\em Hexagons--4 around a point.} (1) Start with a hollow tetrahedron
and cut its corners off.  Cutting the corners off each triangular face
creates 4 hexagonal exterior faces and 4 hexagonal interior faces,
with 4 triangular holes connecting the interior and exterior.  Each
vertex connects 2 exterior and 2 interior hexagons.  (2) Start with a
hollow icosahedron and cut its corners off.  This creates 20 exterior
hexagonal faces, 20 interior hexagonal faces and 12 pentagonal holes
connecting the interior and exterior. (3) Start with a plane
tessellation of hexagons.  It is possible to remove every third
hexagon in such a way that hexagonal holes are created in the plane of
hexagons such that every vertex borders one of the hexagonal holes.
These holes connect to the other side of the plane, so an ant circling
the vertex will visit 2 hexagons on the top of the plane, go through
the hole and traverse 2 hexagons on the bottom of the plane before
returning to where it started.

\item  
{\em Hexagons--8 around a point.}  Tessellate the entire three
dimensional space with octahedrons and tetrahedrons.  This tessellates
space into cells with all triangular faces.  Now cut the corners off
each of the triangular faces turning each of them into regular
hexagons.  Each edge in the original tessellation is truncated to 1/3
its former length as the corners of the triangular faces are cut off.
At each end of one of these truncated original edges is a vertex of
the envelope polyhedron.  Each original edge is bordered by two
octahedral and two tetrahedral volumes, so four back to back
hexagons join at a truncated original edge, and intersect at the
vertex at the end of that edge.  As an ant circles this vertex he will
traverse 2 interior hexagons of an octahedral volume, then 2 interior
hexagons of a tetrahedral volume, then 2 interior hexagons of an
octahedral volume, then 2 interior hexagons of a tetrahedral volume--8
hexagons around a point.  This envelope polyhedron is reminiscent of
{\em octagons--8 around a point} which was based on a tessellation of space
by cubes.

\item
{\em Octagons---4 around a point.} (1) Start with the regular skew
polyhedron {\em squares--6 around a point}, and cut the corners off all the
squares, making all of them into octagons.  Where two squares formerly
met at an edge, two octagons from one side of the surface would meet
at a shortened edge, and each vertex at the end of each of these
shortened edges will now connect these two octagons with two more from
the other side of the surface of the original regular skew
polyhedron. (2) Start with a plane tessellation of squares, truncate
all the squares to create octagons and leave square holes in the
plane. Vertices connect 2 octagons on the top side of the plane with 2
octagons on the bottom side.

\item  
{\em Decagons--4 around a point.} Start with a hollow dodecahedron and
cut off the corners.  This turns each pentagonal face into a decagon
with 10 sides.  The cut off corners become 20 triangular holes,
connecting 12 interior decagons with 12 exterior decagons.  This has
the topology of a sphere with 19 handles, the most complicated
multiply-connected topology of any finite envelope polyhedra.

\item  
{\em Dodecagons--4 around a point.} (1) Start with the regular skew
polyhedron, hexagons--4 around a point.  Cut off the corners of all
the hexagons turning them into dodecagons (with 12 sides each).  Where
two hexagons met on an edge, the edge will now be shortened, and a new
vertex will be created at each end point of each shortened original
edge.  This vertex will now connect two dodecagons from one side of
the original regular skew polyhedron with two dodecagons from the
other side.  That makes {\em dodecagons--4 around a point}. (2) Do the same
operation starting with the regular skew polyhedron {\em hexagons--6 around
a point}. (3) Start with the semi-regular pseudopolyhedron $3^3\times 12^2$ on
page 30 of WBK and remove the triangles. This creates an envelope
polyhedron {\em dodecagons--4 around a point} where a single plane of
decagons is punctuated by triangular holes. (Dodecagons and triangles
tessellate a plane with one triangle and two dodecagons around each
vertex, once the triangles are holes, there are two dodecagons with
top and bottom faces left around each vertex to give 4 around each
vertex).
  
\end{itemize}

\section{Additional Envelope Polyhedra with Mirror Vertices}

Below, in the WBK nomenclature are additional envelope polyhedra
containing mirror image vertices.  Page references are from WBK.

$\mbox{\bf{3}}^{\mbox{\scriptsize{\bf{12}}}}$: (1) From $3^6\times 6$ page 33, by deleting the hexagons.
(2) From $6^6\times 6$
page 34, by deleting the hexagons. (3) Start with Wells's $3^7$ shown on
page 85 of WBK at bottom right. These are icosahedrons connected to
each other by octahedral tunnels in a structure where the icosahedrons
resemble carbon atoms in a diamond structure and the four octahedral
tunnels originating from each resemble the carbon bonds in the diamond
structure. Remove four triangles from each icosahedron, which are
opposite the four octahedral tunnels. This creates a $3^{12}$ envelope
polyhedron with mirror image vertices where there are three exterior
and three interior octahedral tunnel triangles, and three exterior and
three interior icosahedral triangles around each vertex with a
triangular hole connecting the interior and exterior triangles around
each vertex. (4) From $3^6\times4$ on page 95 delete the squares. This leaves
snub cubes missing their square faces linked by open-ended octagonal
tunnels.

$\mbox{\bf{4}}^{\mbox{\scriptsize{\bf{6}}}}$: (1) From $4^4$ page
VIII, 1, by deleting a square from each vertex of an infinitely tall
cylinder with $2n$ sides ($n \geq 2$) to make a series of holes.
There are several ways to do this; see the punctuated planes {\em
  squares--6 around a point} (1) described in the previous section of
the Appendix.  (2) From
$4^3\times 8$ page 3 bottom right, by deleting the octagons. (3) From
$4^3\times 12$ page 10, by deleting the dodecagons. (4) from
$4^3\times 6$ page 11, by deleting the hexagons. (5) From $3^3\times
4^3$ page 13, by deleting the triangles to make a ladder made of
squares with holes in the side rails and cubic rungs with open
ends. This is a vertical stack
of cubes with alternately open east-west ends and open north-south
ends.  (6) From $4^3\times8$ page 19, by deleting the octagons. (7)
From $4^3\times 12$ on page 41, by deleting the dodecagons. (8) From
$3^3\times 4^3$ page 48, by deleting the triangles. (9) From
$3^3\times 4^3$ on page 61 delete the triangles. (10) Tessellate the plane
with dodecagons, squares and hexagons, build dodecagonal prisms, cubes
and hexagonal prisms on these.  Delete their tops and bottoms.  (11) From $4^3\times 6$ on page 72 delete the hexagons. (12)
From $4^3\times 6$ on page 86 delete the hexagons. (13) From
$4^3\times 8$ on page 89 delete the octagons. (14) From $4^3\times 8$
on page 98 delete the octagons.  (15) From $(4\times 8)^2$ on page 97,
deleted the octagons and add a square at each vertex connecting the
remaining two squares to make a bent \textsf{L} shaped pattern of 3
squares (where the two legs of the \textsf{L} are bent in opposite
directions at the edges of the squares by $45^\circ$ each).  This
gets traversed twice by an ant circling the vertex to give $4^8$.  The
structure consists of octagonal prisms, missing their octagonal tops
and bottoms, pasted together at right angles at square
sides---octagonal rings of squares meeting at right angles in a
three-dimensional structure.  (This can also be constructed from {\em
  squares--10 around a point} mentioned in the main body of this paper
by deletion of the appropriate squares).  

$\mbox{\bf{4}}^{\mbox{\bf{\scriptsize{8}}}}$: (1) From $4^5$ page 17,
by deleting squares. Bottom right picture shows some squares face on;
these represent towers of open cubes seen from the top, eliminate
squares from these towers seen edge on in the bottom right picture, so
as to leave the squares seen face on in the bottom right picture as
fins connecting the other 2 by 2 boxlike structures. (2, 3, 4, \& 5)
From $4^5$ on page 20, delete the horizontal squares in the right hand
picture. Or from $4^5$ on page 20, instead delete alternate squares in
the filmstrips of squares connecting octagonal columns of
squares. There are two ways to do this. Or delete the squares on the
vertical sides of half of the open cubes connecting the octagonal
columns. (6) From $4^4\times 6$ on page 37, by deleting the
hexagons. (7) From $4^4\times 6$ on page 43, by deleting the
hexagons. (8, 9, 10, \& 11) From $4^5$ on page 46, delete the squares
seen face-on in the upper right hand picture. Or from $4^5$ on page 46
instead delete alternate squares in the filmstrips of squares
connecting hexagonal and octagonal columns of squares. There are two
ways to do this. Or delete squares on the vertical sides of the open
cubes connecting the hexagonal columns in the bottom right
figure. (12) From $4^4\times 8$ on page 49, by deleting the
octagons. (13) From $4^4\times 8$ on page 55 delete the octagons. (14)
From $4^4\times 8$ on page 56 delete the octagons. (15) From
$4^4\times 6$ on page 57 delete the hexagons. (16) From $4^4\times$ 6
on page 58 delete the hexagons. (17) From $4^4\times 12$ on page 59
delete the dodecagons. (18) From $3^3\times 4^4$ on page 60 at left
delete the triangles, leaving two parallel layers of checkerboard
pattern of just the white squares connected by filmstrips of
squares. (19) From $3^3\times 4^4$ on page 62 delete the triangles.
(20, 21, 22) Attach square fins (one per vertex) to the single layer
of $4^6$ (10) described above to attach to layers above to make a
vertical stack of layers.  There are three ways to do this, as three
double-sided squares meet at a vertex in a single layer, and so we
have three different places to attach a fin.  (23) Octagons and
squares tile the plane.  Turn them into prisms and make a planar layer
of octagonal prisms and cubes without their tops and bottoms.  Add a
series of parallel square fins above opposite sides of each cube.  Repeat
so that these parallel fins connect to a similar layer above with
mirror image vertices.

$\mbox{\bf{4}}^{\mbox{\scriptsize{\bf{10}}}}$:
(1) From $4^6$ on page 45, delete every other square forming the
hexagonal rings. (2) From $4^6$ on page 53, delete horizontal squares
from North-South ties in the lower left figure. (3) Delete similar
squares from $4^6$ on page 54 at upper right. (4) Delete horizontal squares on
North-South tunnels from Petrie $4^6$ on page 67, upper right
figure.  (5) Take {\em squares--8 around a point} (8) described
above and add fins connecting it to other copies of {\em squares--8
  around a point}. This can also be obtained from Fig.~5 by deleting
the 12 squares at the back of Fig.~5, leaving one row of open ended
triangular prisms with fins sticking up. This structure repeats
vertically. It may also be viewed as two vertical checkerboards
punctuated with holes [{\em squares--6 around a point} (1) described
  above] connected by horizontal east-west zig-zag filmstrips [{\em
    squares--4 around a point}] giving $4^{10}$.  The holes in the {\em
  squares--6 around a point} punctured checkerboards jog back and
forth half a square from layer to layer. (6) Start with {\em
  squares--10 around a point} (2). This has vertical checkerboards
punctuated by square holes, connected by horizontal sequences of
squares which are sides of cubes. Now between checkerboards rotate by
$90^\circ$ to produce vertical sequences of squares that are sides of
cubes. We now have punctuated checkerboards connected by alternating
horizontal and vertical sequences of squares that are sides of cubes.

$\mbox{\bf{4}}^{\mbox{\scriptsize{\bf{12}}}}$: (1 \& 2) Start with
{\em squares--10 around a point} (3) and add fins to the top of the
single layer of open ended triangular and hexagonal prisms tiling a
plane. There are two ways to do this. Add vertical fins to alternate
vertical square sides of the hexagonal prisms, or add vertical fins to
square sides between two triangle prisms where the fin is in the same
plane as sides of the two hexagonal prisms that it touches at
vertices. Now repeat vertically to create alternate layers of open
prisms and fins. These of course have mirror vertices. (3, 4) Start
with the plane tessellation $3^3 \times 4^2$ and construct prisms over
it. Now take the tops and bottoms off these prisms. These are cubes
and triangular prisms. Add square fins to this to connect to the next
vertical layer of open prisms. There are two ways to do this that
produce mirror vertices. Place the fins in parallel above alternate
squares joining a cube and a triangular prism, or place the fins
instead above squares connecting triangular prisms

$\mbox{\bf{4}}^{\mbox{\scriptsize{\bf{14}}}}$: In Fig.~5 the repeated
fins and rows of squares repeated vertically produce vertical
checkerboards punctuated with holes [{\em squares--6 around a point}
  (1)] connected by horizontal zig-zag filmstrips of squares. Now
rotate the set of 12 squares at the back of the figure by $90^\circ$
in a vertical plane, so that the zig-zag filmstrip of squares stands
vertically. We now have a series of punctuated checkerboards connected
to adjacent checkerboards by horizontal rows of zig-zag filmstrips of
{\em squares--4 around a point} on one side and vertical rows of
zig-zag filmstrips of {\em squares--4 around a point} on the other
side. This gives 4 + 6 + 4 = 14 squares around a point: $4^{14}$.

Interestingly, squares, which have 4 sides, seem
to produce the greatest variety of structures in 3D, just as carbon,
which has four bonds with adjacent atoms, produces the richest
chemistry.

\end{document}